\pdfoutput=1
\documentclass[letterpaper, 10 pt, conference]{ieeeconf}

\IEEEoverridecommandlockouts
\usepackage{amsmath,amssymb,amsfonts,amsthm}
\usepackage{graphicx}
\usepackage{cite}
\usepackage{xcolor}
\usepackage{tikz}
\usetikzlibrary{arrows.meta,positioning,calc,fit,decorations.pathreplacing}
\tikzset{
  agent/.style={circle,draw,thick,minimum size=5.5mm,inner sep=0pt,font=\scriptsize},
  inn/.style={agent,fill=blue!15},
  outn/.style={agent,fill=red!15},
  flatn/.style={agent,fill=green!20},
  blk/.style={draw,thick,minimum height=6mm,minimum width=7mm,inner sep=2pt},
  sumn/.style={circle,draw,thick,minimum size=3.5mm,inner sep=0pt},
  sig/.style={-{Latex[length=1.6mm]},thick},
}

\newtheorem{proposition}{Proposition}
\newtheorem{lemma}{Lemma}
\newtheorem{corollary}{Corollary}
\theoremstyle{definition}
\newtheorem{definition}{Definition}
\newtheorem{assumption}{Assumption}
\newtheorem{problem}{Problem}
\newtheorem{examplex}{Example}
\newenvironment{example}{\pushQED{\qed}\examplex}{\popQED\endexamplex}
\theoremstyle{remark}
\newtheorem{remarkx}{Remark}
\newenvironment{remark}{\pushQED{\qed}\remarkx}{\popQED\endremarkx}

\newcommand{\R}{\mathbb{R}}
\newcommand{\G}{\mathcal{G}}
\newcommand{\V}{\mathcal{V}}
\newcommand{\E}{\mathcal{E}}
\newcommand{\N}{\mathcal{N}}
\newcommand{\one}{\mathbf{1}}
\newcommand{\e}{\mathrm{e}}
\newcommand{\dist}{\mathrm{dist}}
\newcommand{\diag}{\mathrm{diag}}

\newcommand{\Sg}{\Sigma_{\G}}

\title{\LARGE \bf Differential Flatness of Nonlinear Consensus Networks: \\A Graph-Theoretic Perspective}

\author{Daniel Zelazo$^{1}$ and Timm Faulwasser$^{2}$%
\thanks{The work of Daniel Zelazo was supported by the Israel Science Foundation (ISF) under grant no.\,453/24. The work of Timm Faulwasser was supported by the Deutsche Forschungsgemeinschaft (DFG) under project no.\,527447339.}%
\thanks{$^{1}$D. Zelazo is with the Stephen B. Klein Faculty of Aerospace Engineering, Technion -- Israel Institute of Technology, Haifa, Israel. {\tt\small dzelazo@technion.ac.il}}%
\thanks{$^{2}$T. Faulwasser is with the Institute of Control Systems, Hamburg University of Technology, Hamburg, Germany. {\tt\small timm.faulwasser@ieee.org}}%
}

\begin{document}
\maketitle
\thispagestyle{empty}
\pagestyle{empty}

\begin{abstract}
We study the property of differential flatness for networks of identical flat agents coupled through a consensus protocol, where a single exogenous input acts on one agent and the output of another agent is measured. We distinguish whether the network admits a flat output at all from whether the measured task output is a flat output. For the former, we show that flatness at a synchronized equilibrium requires controllability of the Laplacian leader--follower pair, which is also sufficient for linear agents but not for nonlinear ones. For the latter, we show that the relative degree of the task output is determined by its graph distance to the input node, so that it is a flat output if and only if the graph is a path with the input and output nodes at its ends. Numerical examples illustrate how flatness enables finite-time transitions between consensus configurations through a single input.
\end{abstract}

\section{Introduction}
\label{sec:intro}

In the canonical \emph{consensus} protocol~\cite{Mesbahi2010}, each agent drives its output towards those of its neighbors through diffusive coupling, and the resulting collective behavior is dictated by the properties of the graph Laplacian. In many applications, a network may be required to track an external reference known only to a subset of the agents~\cite{Ren2007,Cao2012}, or reconfigure their spatial formation~\cite{Turpin2012,Oh2015}. Often such interaction with the network is only possible through a few \emph{leader} nodes~\cite{Rahmani2009,Egerstedt2012}, and only a few nodes can be measured. This raises a structural question: which trajectories of a network output can be generated by an input applied at a remote node?

For a single system, \emph{differential flatness}~\cite{Fliess1995} provides a powerful answer to this question. If a system is flat, there exists an output, the \emph{flat output}, such that all states and inputs can be expressed algebraically as functions of this output and finitely many of its time derivatives. Trajectory generation then reduces to planning a sufficiently smooth curve for the flat output. This property is the foundation for a large body of work on motion planning,  feedforward control, constrained trajectory generation and reachability~\cite{Martin2003,Levine2009,Faulwasser2014}. It is also the core mechanism behind computed-torque concepts in robotics. For linear systems, flatness is equivalent to controllability~\cite{Rothfuss97b,Levine2009}, a flat output being furnished by the Brunovsk\'y canonical form~\cite{Martin2003,Brunovsky1970}. For SISO nonlinear systems, flatness is equivalent to static feedback linearizability~\cite{Fliess1995,Charlet1989,Jakubczyk1980}, and flat outputs are exactly the outputs of full relative degree~\cite{Isidori1995}; in particular, they are essentially unique up to a reparametrization.

Consensus is well understood for single and higher-order integrators and passive agents~\cite{Mesbahi2010,Ren2007,Chopra2006}, but these works concern the stability of agreement rather than steering the network through a few input nodes. Controllability of leader--follower consensus networks was introduced in~\cite{Tanner2004}; symmetries of the graph with respect to the leader render it uncontrollable~\cite{Rahmani2009,Egerstedt2012}, and distance partitions bound the controllable subspace~\cite{Zhang2014}. For linear networks with a single input and a single output node, the relative degree is determined by the graph distance~\cite{Briegel2011}, and the transfer zeros are related to paths in the graph~\cite{AbadTorres2015}. Flatness, in turn, is widely used for multi-agent coordination, e.g., in formation flight~\cite{Turpin2012}, but there each agent has its own input and is flat individually; similarly, \cite{Yang2026} identifies couplings between flat subsystems, each with its own input, that preserve flatness of the joint system. In contrast, we study networks of flat agents diffusively coupled through the flat outputs and driven by a \emph{single} network input.  We investigate \emph{flatness} of such a network as \emph{an emergent property induced by the interconnection structure}. To the best of our knowledge, this interplay between graph structure and differential flatness has not been addressed.

Specifically, we ask when a \emph{network of nonlinear SISO flat agents}, coupled through a consensus protocol and driven from a single node, is itself differentially flat, and how this depends on where the network is actuated and where it is observed. Given the agents and the interaction graph, two design questions arise:
\begin{itemize}
    \item[(i)] \emph{Actuator placement}: at which node can a single input be injected so that the network is flat?
    \item[(ii)] \emph{Sensor placement}: where must the \emph{task output}  be located, relative to a given input node, for it to be a flat output?
\end{itemize}

To tackle both questions, we distinguish the property of \emph{network flatness}, i.e., the existence of a flat output, which constrains where the input can be placed, from \emph{flatness of the task output}.
We show that the two guiding questions have different answers: Actuator placement is constrained by the controllability of Laplacian leader--follower dynamics, which is necessary for flatness and sufficient for linear agents, so that graph symmetries about the input node rule out flatness. Classical symmetry obstructions to leader--follower controllability~\cite{Rahmani2009} are therefore obstructions to flatness. We show by a counterexample that this condition is not sufficient for nonlinear systems.

The task output, by contrast, is a flat output only if its graph distance from the input node is one less than the number of agents, which forces the graph to be a path with the input and the task output at its two ends. Put differently, we prove that the latter property is much stronger.
For the task output, we show that the relative degree from the input node to the output node is determined by their graph distance, increasing by the agent order with each edge, which generalizes the linear first-order result of~\cite{Briegel2011}. In every other flat network, the task output is not a flat output, and the flat output may even be a combination of several agents' outputs, as we illustrate with an explicit example.

The remainder of the paper is organized as follows. Sections~\ref{sec:prelim}--\ref{sec:problem} give background and the problem statement, Section~\ref{sec:main} graph conditions for both flatness aspects of diffusively coupled networks of flat agents, and Section~\ref{sec:example} numerical examples of trajectory planning.

\paragraph*{Notation} Throughout, $\e_i$ denotes the $i$-th unit coordinate vector and $I_n$ the $n\times n$ identity matrix, with the dimension of $\e_i$ clear from context, $\one_N\in\R^N$ is the vector of ones, and $\otimes$ is the Kronecker product. We use standard graph-theoretic notation~\cite{Mesbahi2010}: $\G=(\V,\E,w)$ is an undirected graph with nodes $\V=\{1,\dots,N\}$, edges $\E$ and weights $w_{ij}=w_{ji}>0$, $\N_i$ is the neighborhood of node $i$, and $L\in\R^{N\times N}$ is the weighted Laplacian, which is symmetric positive semidefinite with $L\one_N=0$. The distance $\dist(i,j)$ is the number of edges of a shortest path between nodes $i$ and $j$, and $\mathcal P_N$ is the path graph on $N$ nodes. For a smooth function $y(t)$, $y^{(k)}$ denotes its $k$-th time derivative, and $\bar y^{[k]} := (y,\dot y,\dots,y^{(k)})$. For an open set $U\subseteq\R^n$, smooth vector fields $f,g:U\to\R^n$ and a smooth function $h:U\to\R$, the differential of $h$ is the row vector $\mathrm dh=\partial h/\partial x$, the Lie derivative of $h$ along $f$ is $\mathcal L_fh=\mathrm d h\,f$, with $\mathcal L_f^0h=h$ and $\mathcal L_f^{k+1}h=\mathcal L_f(\mathcal L_f^kh)$, and the Lie bracket of $f$ and $g$ is $[f,g]=\frac{\partial g}{\partial x}f-\frac{\partial f}{\partial x}g$, with $\mathrm{ad}_f^0g=g$ and $\mathrm{ad}_f^{k+1}g=[f,\mathrm{ad}_f^kg]$. Smooth vector fields $v_1,\dots,v_k$ on $U$ define the distribution $\Delta(x)=\mathrm{span}\{v_1(x),\dots,v_k(x)\}\subseteq\R^n$, $x\in U$, which is called nonsingular if $\dim\Delta(x)$ is constant on $U$, and involutive if $[v_i,v_j](x)\in\Delta(x)$ for all $i,j\in\{1,\dots,k\}$ and all $x\in U$. We write $\mathrm dh\perp\Delta$ if $\mathrm dh(x)\,v=0$ for all $v\in\Delta(x)$ and all $x\in U$.
\section{Preliminaries}
\label{sec:prelim}

We first recall the consensus protocol and the controllability of leader--follower consensus dynamics, which will govern the existence of a flat output for the network. We then review differential flatness, with emphasis on single-input systems, for which existence of a flat output and flatness of a given output are characterized by the geometric controllability distributions and by relative degree, respectively.

\subsection{Leader--follower consensus}

For single-integrator agents $\dot y_i=u_i$, the consensus protocol
\begin{equation}\label{eq:consensus}
u_i = \sum_{j\in\N_i} w_{ij}\,(y_j-y_i), \quad\text{i.e.,}\qquad \dot y = -L y,
\end{equation}
drives all outputs to the average of their initial values whenever $\G$ is connected~\cite{OlfatiSaber2004}. When an exogenous input $u_G$ is injected at a leader node $\ell$, one obtains the \emph{leader--follower} dynamics
\begin{equation}\label{eq:leaderfollower}
    \dot y=-Ly+\e_\ell u_G.
\end{equation}
Its controllability is characterized by the spectrum of $L$ and will reappear in Section~\ref{sec:main} as a necessary condition for network flatness.
\begin{lemma}[Leader--follower controllability~\cite{Rahmani2009,Egerstedt2012}]\label{lem:leaderfollower}
   System~\eqref{eq:leaderfollower} is controllable if and only if no eigenvector of $L$ is orthogonal to $\e_\ell$ and all eigenvalues of $L$ are simple.
\end{lemma}

\subsection{Differential flatness}

Consider the nonlinear system
\begin{equation}\label{eq:gen_sys}
\dot x = f(x,u),\quad x\in\R^n,\; u\in\R^p,
\end{equation}
with $f:\R^n\times \R^p \to \R^n$ a smooth vector field.

\begin{definition}[Differential flatness~\cite{Fliess1995}]\label{def:flat}
System~\eqref{eq:gen_sys} is \emph{(differentially) flat} if there exist an output $z\in\R^p$, integers $r\ge1$, $s\ge 0$, and smooth maps $h,\Phi,\Psi$ such that, along all trajectories,
\begin{subequations}\label{eq:flat_maps}
\begin{align}
z &= h\big(x,u,\dot u,\dots,u^{(s)}\big),\\
x &= \Phi\big(z,\dot z,\dots,z^{(r-1)}\big)=: \Phi\left(\bar z^{[r-1]}\right),\\
u &= \Psi\big(z,\dot z,\dots,z^{(r)}\big)=: \Psi\left(\bar z^{[r]}\right).
\end{align}
\end{subequations}
holds and the components of $z$ do not satisfy any differential equation.
The output $z$ is then called a \emph{flat output}.
\end{definition}

For smooth $z_{\rm ref}$, $x_{\rm ref}=\Phi(\bar z_{\rm ref}^{[r-1]})$, $u_{\rm ref}=\Psi(\bar z_{\rm ref}^{[r]})$ is a trajectory of~\eqref{eq:gen_sys}, so trajectory generation can be simplified, see, e.g.,~\cite{Martin2003,Faulwasser2014}.

\begin{lemma}[Linear systems~\cite{Fliess1995,Brunovsky1970}]\label{lem:lin_flat}
A linear system $\dot x=Ax+Bu$ is flat if and only if $(A,B)$ is controllable. For $p=1$, an output $z=\gamma^\top x$ is a flat output if and only if $\gamma^\top=\kappa\,\e_n^\top\mathbf C^{-1}$ for some $\kappa\neq0$, where $\mathbf C=[B,AB,\dots,A^{n-1}B]$; i.e., the linear flat output is unique up to scaling.
\end{lemma}

Let \eqref{eq:gen_sys} be single-input control-affine
\begin{equation}\label{eq:affine}
\dot x = f(x)+g(x)u,\quad \dim u=1.
\end{equation}
Recall that an output $y=h(x), \dim y=1$ has \emph{relative degree} $\rho$ at $x^0$ if $\mathcal L_g\mathcal L_f^kh(x)=0$ for all $k<\rho-1$ and all $x$ near $x^0$, and $\mathcal L_g\mathcal L_f^{\rho-1}h(x^0)\neq 0$~\cite{Isidori1995}. Then $y^{(k)}=\mathcal L_f^kh(x)$ for $k<\rho$ and $y^{(\rho)}=\mathcal L_f^\rho h(x)+\mathcal L_g\mathcal L_f^{\rho-1}h(x)\,u$.

\begin{lemma}[Nonlinear single-input systems \cite{Fliess1995,Martin2003,Levine2009}]\label{lem:siso_flat}
Consider~\eqref{eq:affine} around $x^0$. Then:
\begin{enumerate}
  \item[(a)] the system is locally flat around $x^0$ if and only if it is locally static feedback linearizable, i.e., if and only if the distributions $\Delta_k=\mathrm{span}\{g,\mathrm{ad}_fg,\dots,\mathrm{ad}_f^{k}g\}$ are involutive near $x^0$ for $k=0,\dots,n-2$ and $\dim\Delta_{n-1}(x^0)=n$;
  \item[(b)] the output $y=h(x)$ is a local flat output if and only if it has relative degree $\rho=n$ at $x^0$, in which case $x\mapsto(h,\mathcal L_fh,\dots,\mathcal L_f^{n-1}h)$ is a local diffeomorphism;
  \item[(c)] if $h_1(x)$ and $h_2(x)$ are both local flat outputs, then $h_2=\theta\circ h_1$ near $x^0$ for some smooth $\theta$ with $\theta'\neq0$.
\end{enumerate}
\end{lemma}

The two conditions in~(a) play different roles. The rank condition on $\Delta_{n-1}$ is a controllability condition, whereas involutivity is what makes a flat output exist. Under the rank condition, $\Delta_{n-2}$ is nonsingular of dimension $n-1$ near $x^0$, and by Frobenius' theorem~\cite{Isidori1995} it is involutive near $x^0$ if and only if there is a smooth function $h$ with $\mathrm dh(x^0)\neq0$ and $\mathrm dh\perp\Delta_{n-2}$ near $x^0$, which, again under the rank condition, is equivalent to relative degree $n$ at $x^0$. If $\Delta_{n-2}$ is not involutive, any $h$ with $\mathrm dh\perp\Delta_{n-2}$ must also annihilate the brackets of its vector fields, so that $\mathrm dh$ vanishes wherever these brackets together with $\Delta_{n-2}$ span $\R^n$. Hence, unlike for linear systems, controllability is necessary but not sufficient for flatness. Part~(c) follows from~(b), since the annihilator of $\Delta_{n-2}$ is one-dimensional.

\begin{remark}[Existence versus identification]\label{rem:lin_nec}
Lemma~\ref{lem:siso_flat} separates two questions. Part~(a) concerns the \emph{existence} of a flat output and depends on $(f,g)$ only, whereas part~(b) concerns whether a \emph{given} output $h$ is flat, which by~(c) holds if and only if $h=\theta\circ h^\ast$ for an essentially unique flat output $y^\ast=h^\ast(x)$. If $x^0$ is an equilibrium for $u=0$, then $\mathrm{ad}_f^kg(x^0)=(-1)^kA^kb_0$ with $A=\partial f/\partial x(x^0)$ and $b_0=g(x^0)$, so the rank condition in~(a) is controllability of the linearization, which is therefore necessary for local flatness. For the network, the two questions become the two problems introduced in Section~\ref{sec:problem}.
\end{remark}

\section{Problem Statement}
\label{sec:problem}

The building blocks of the studied network are agents that are differentially flat in isolation, with their output as flat output.
Specifically, we consider $N$ identical SISO agents in nonlinear normal form,
\begin{equation}\label{eq:agent}
\begin{aligned}
\dot x_{i,k} &= x_{i,k+1}, \qquad k=1,\dots,q-1,\\
\dot x_{i,q} &= a(x_i)+b(x_i)\,u_i,\\
y_i &= x_{i,1},
\end{aligned}
\end{equation}
with state $x_i=(x_{i,1},\dots,x_{i,q})\in\R^q$, input $u_i\in\R$ and output $y_i\in\R$, $i\in\V$. In compact form, the agent dynamics read
\begin{equation}\label{eq:agent_compact}
\dot x_i=A_cx_i+\e_q\big(a(x_i)+b(x_i)\,u_i\big), \quad y_i=\e_1^\top x_i,
\end{equation}
where $A_c=\left[\begin{smallmatrix}0&I_{q-1}\\0&0\end{smallmatrix}\right]\in\R^{q\times q}$ and $\e_1,\e_q\in\R^q$.
\begin{assumption}\label{ass:agent}
The functions $a,b:\R^q\to\R$ are smooth and $b(x_i)>0$ for all $x_i\in\R^q$.
\end{assumption}

Under Assumption~\ref{ass:agent}, each agent is (globally) flat with flat output $y_i$: $x_i=\bar y_i^{[q-1]}$ and $u_i=(y_i^{(q)}-a(\bar y_i^{[q-1]}))/b(\bar y_i^{[q-1]})$.

\begin{remark}[Generality of~\eqref{eq:agent}]\label{rem:normal_form}
By Lemma~\ref{lem:siso_flat}(b), a SISO control-affine agent with flat output $ h$ takes the form~\eqref{eq:agent}, at least locally, in the coordinates $\tilde x=(h,\mathcal L_{f} h,\dots,\mathcal L_{ f}^{q-1} h)$. Hence~\eqref{eq:agent} covers any network of SISO flat agents that exchange information via their \emph{flat outputs}.
\end{remark}

\subsection{Consensus coupling and the network system}

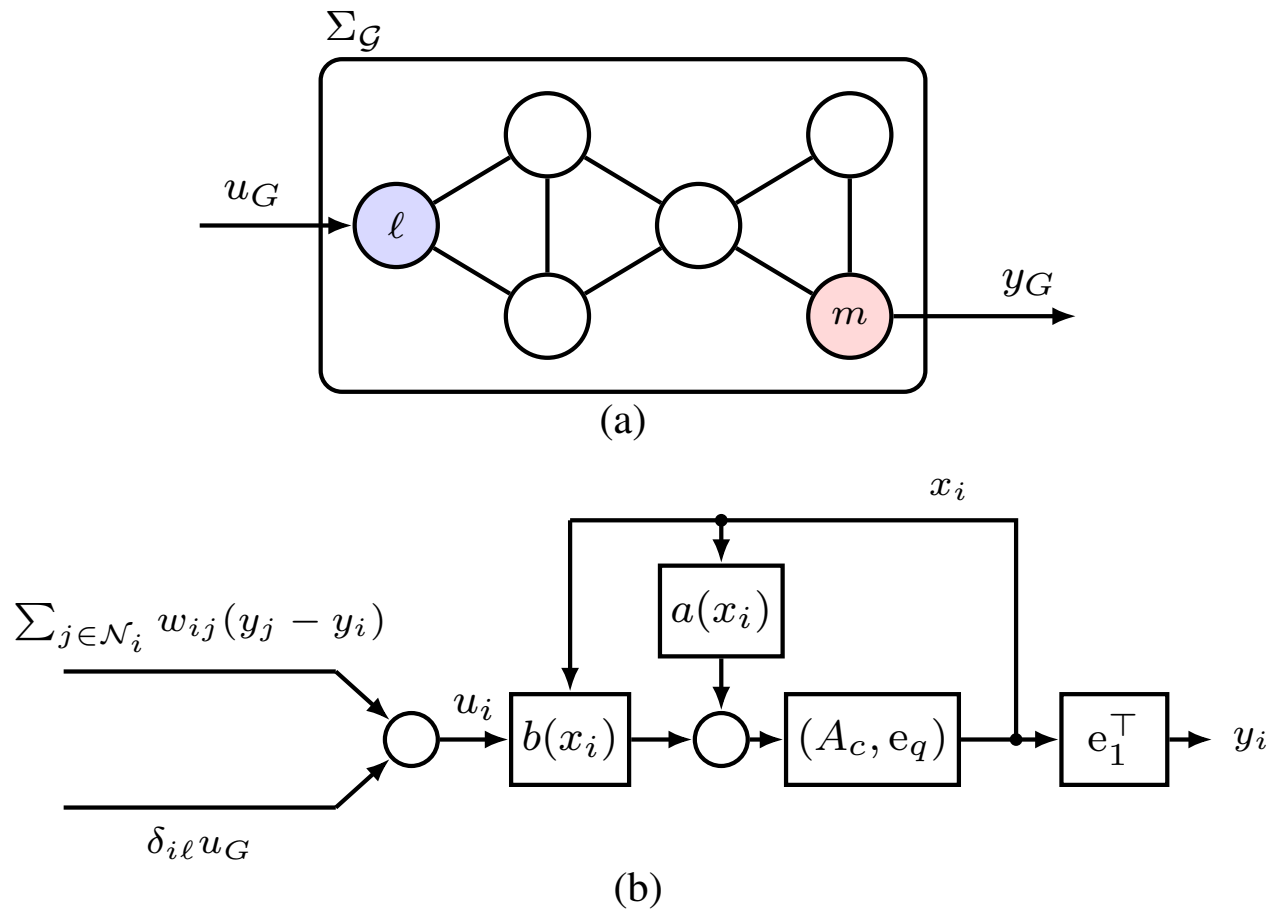
\begin{figure}[t]
\begin{center}
\begin{tikzpicture}[font=\footnotesize]
\node[inn]   (n1) at (0,0)      {$\ell$};
\node[agent] (n2) at (1.0,0.6)  {};
\node[agent] (n3) at (1.0,-0.6) {};
\node[agent] (n4) at (2.0,0)    {};
\node[agent] (n5) at (3.0,0.6)  {};
\node[outn]  (n6) at (3.0,-0.6) {$m$};
\draw[thick] (n1)--(n2) (n1)--(n3) (n2)--(n3) (n2)--(n4) (n3)--(n4)
             (n4)--(n5) (n4)--(n6) (n5)--(n6);
\node[draw,rounded corners,thick,fit=(n1)(n2)(n3)(n4)(n5)(n6),inner sep=6pt] (box) {};
\node[anchor=south west,inner sep=1pt] at (box.north west) {$\Sigma_{\mathcal G}$};
\draw[sig] (-1.3,0) -- node[above,pos=0.35]{$u_G$} (n1.west);
\draw[sig] (n6.east) -- (4.5,-0.6) node[above left]{$y_G$};
\node at (1.5,-1.3) {(a)};
\begin{scope}[shift={(0.1,-3.4)}]
\node[sumn] (s1) at (0,0) {};
\node[blk]  (bb) at (1.05,0) {$b(x_i)$};
\node[sumn] (s2) at (2.05,0) {};
\node[blk]  (ss) at (3.05,0) {$(A_c,\e_q)$};
\node[blk]  (cc) at (4.65,0) {$\e_1^\top$};
\node[blk]  (aa) at (2.05,0.85) {$a(x_i)$};
\coordinate (tap) at (4.0,0);
\draw[thick] (-2.3,0.45) -- node[above,font=\scriptsize]{$\sum_{j\in\mathcal N_i}w_{ij}(y_j-y_i)$} (-0.5,0.45);
\draw[sig] (-0.5,0.45) -- (s1);
\draw[thick] (-2.3,-0.45) -- node[below,font=\scriptsize]{$\delta_{i\ell}u_G$} (-0.5,-0.45);
\draw[sig] (-0.5,-0.45) -- (s1);
\draw[sig] (s1) -- node[above]{$u_i$} (bb);
\draw[sig] (bb) -- (s2);
\draw[sig] (aa) -- (s2);
\draw[sig] (s2) -- (ss);
\draw[sig] (ss) -- (cc);
\draw[sig] (cc) -- (5.3,0) node[right,font=\scriptsize]{$y_i$};
\coordinate (fb) at (0,1.45);
\fill (tap) circle (1.2pt);
\draw[thick] (tap) -- (tap |- fb) -- node[above,pos=0.15,font=\scriptsize]{$x_i$} (bb.north |- fb);
\draw[sig] (bb.north |- fb) -- (bb.north);
\draw[sig] (aa.north |- fb) -- (aa.north);
\fill (aa.north |- fb) circle (1.2pt);
\node at (1.5,-1.0) {(b)};
\end{scope}
\end{tikzpicture}
\end{center}
\caption{(a) Network $\Sg$: agents coupled by the consensus protocol~\eqref{eq:coupling}; the exogenous input $u_G$ acts at node $\ell$, and the output $y_G=y_m$ of node $m$ is measured. (b) Agent $i$ in the normal form~\eqref{eq:agent_compact}, with its flat output $y_i=x_{i,1}$ exchanged with its neighbors.}
\label{fig:setup}
\end{figure}

The agents interact over a connected undirected graph $\G$ through the consensus protocol~\eqref{eq:consensus} applied to their flat outputs. In addition, an exogenous input $u_G\in\R$ acts on a single \emph{input node} $\ell\in\V$, and the output $y_G\in\R$ of a single \emph{output node} $m\in\V$ is measured:
\begin{equation}\label{eq:coupling}
u_i = \sum_{j\in\N_i}w_{ij}\,(y_j-y_i)+\delta_{i\ell}\,u_G,\quad y_G = y_m,
\end{equation}
where $\delta_{i\ell}$ is the Kronecker delta. Stacking $x=(x_1,\dots,x_N)\in\R^{Nq}$, $\mathbf a(x)=(a(x_1),\dots,a(x_N))$, and $B(x)=\diag(b(x_1),\dots,b(x_N))$, the interconnection~\eqref{eq:agent_compact}--\eqref{eq:coupling} reads
\begin{equation}\label{eq:network}
\Sg:\;
\left\{
\begin{aligned}
\dot x &= (I_N\otimes A_c)x + b(x_\ell)\,(\e_\ell\otimes \e_q)\,u_G\\
&\quad + (I_N\otimes \e_q)\big(\mathbf a(x) - B(x)\,L\,(I_N\otimes \e_1^\top)x\big),\\
y_G &= (\e_m^\top\otimes \e_1^\top)\,x,
\end{aligned}
\right.
\end{equation}
Here $\e_\ell,\e_m\in\R^N$ select the input and output agents, whereas $\e_1,\e_q\in\R^q$ select the flat-output and input coordinates within an agent, and $b(x_\ell)$ is the input gain of agent~$\ell$. The setup is depicted in Fig.~\ref{fig:setup}.

For the subsequent analysis we will rewrite \eqref{eq:network} in more compact form as
\begin{subequations}\label{eq:network_affine}
\begin{align}
    \dot x &= \underbrace{\begin{bmatrix}
        F_1(x_1, x_{j \in \mathcal N_1}) \\ F_2(x_2, x_{j \in \mathcal N_2}) \\\vdots \\ F_N(x_N, x_{j \in \mathcal N_N})
    \end{bmatrix}}_{=F(x)} + \underbrace{\begin{bmatrix}
       0\\\vdots\\ b(x_\ell)\\0
    \end{bmatrix}}_{G(x) = b(x_\ell)\e_\ell\otimes \e_q}u_G\\
    y_G &= h(x) = (\e_m^\top\otimes \e_1^\top)\,x = \e_1^\top x_m
\end{align}
\end{subequations}
which is of the form~\eqref{eq:affine} with $f=F$, $g=G$ and $n=Nq$.

\begin{assumption}\label{ass:graph}
The graph $\G$ is connected and $N\ge 2$.
\end{assumption}

For $c\in\R$, the \emph{synchronized configuration} $x^\ast_c:=\one_N\otimes c\,\e_1$ is the state in which all agents agree on the flat output $y_i=c$ and all higher derivatives vanish. Since $(I_N\otimes A_c)x^\ast_c=0$ and $L(I_N\otimes \e_1^\top)x^\ast_c=0$, the configuration $x^\ast_c$ is an equilibrium of~\eqref{eq:network} with $u_G=0$ if and only if $c$ belongs to the set
\begin{equation}\label{eq:sync_set}
\mathbb S:=\{c\in\R:\ a(c\,\e_1)=0\},
\end{equation}
in which case we call $x^\ast_c$ a \emph{synchronized equilibrium}.

\begin{remark}\label{rem:not_stable}
For $q=1$, $a\equiv 0$ and $b\equiv 1$, \eqref{eq:network} reduces to the classical leader--follower consensus network $\dot y=-Ly+\e_\ell u_G$~\cite{Rahmani2009,Egerstedt2012}. For $q\ge 2$, the coupling~\eqref{eq:coupling} need not render agreement attractive; we do not require this, since~\eqref{eq:network} is the object to be \emph{steered} via $u_G$.
\end{remark}

\subsection{Problem formulation}

Each agent is flat, but it is not clear whether the interconnection~\eqref{eq:network}, which has a single input for $Nq$ states, is flat. Following Remark~\ref{rem:lin_nec}, we distinguish two properties of \eqref{eq:network}.

\begin{definition}[Network flatness]\label{def:flat_network} The interconnected system~\eqref{eq:network} is said to be a \emph{differentially flat network} if $\Sg$ admits a flat output $z$ in the sense of Definition~\ref{def:flat}, measured or not.
\end{definition}
 This property depends on $a$, $b$, $\G$ and $\ell$, but not on $m$, and hence is a question of where to place the input. Network flatness suffices for \emph{motion planning}: for every sufficiently smooth $z_{\rm ref}$, the input $u_G=\Psi(\bar z_{\rm ref}^{[Nq]})$ steers the whole network along $x_{\rm ref}=\Phi(\bar z_{\rm ref}^{[Nq-1]})$, e.g., between two consensus configurations.
\begin{definition}[Flat task output]\label{def:flat_output} The interconnected system~\eqref{eq:network} is said to have a \emph{flat task output} if $y_G$ is a flat output.
\end{definition}
The flat output $z$ of Definition~\ref{def:flat_network} is in general a combination of the states of many agents and is neither measured nor meaningful for the task, whereas node $m$ is typically dictated by the application. By Lemma~\ref{lem:siso_flat}(c), task flatness asks whether $y_m$ coincides with the essentially unique flat output of $\Sg$. In addition to planning, a flat task output enables
\begin{itemize}
  \item \emph{algebraic reconstruction} of the entire network state $x=\Phi(\bar y_G^{[Nq-1]})$ from derivatives of a single measured signal, without an observer; and
  \item exact \emph{output tracking}: for every sufficiently smooth $y_{\rm ref}$, $u_G=\Psi(\bar y_{\rm ref}^{[Nq]})$ drives the measured agent exactly along $y_{\rm ref}$, with no internal dynamics left over.
\end{itemize}
Definition~\ref{def:flat_output} implies Definition~\ref{def:flat_network}. This motivates the following problems.

\begin{problem}[Network flatness]\label{prob:network}
Determine graph-theoretic conditions under which $\Sg$ is (locally) flat, i.e., admits \emph{some} flat output, possibly not a node output.
\end{problem}

\begin{problem}[Task flat outputs]\label{prob:node}
Given $\G$, characterize all pairs $(\ell,m)\in\V\times\V$ for which $y_G=y_m$ is a (local) flat output of $\Sg$.
\end{problem}

Both problems are addressed in Section~\ref{sec:main}, and Section~\ref{sec:example} illustrates how the results can be used for trajectory planning.

\section{When Is a Consensus Network Flat?}
\label{sec:main}

We now address Problems~\ref{prob:network} and~\ref{prob:node} in turn, using a linearization argument for the former and a relative degree computation for the latter, and then relate the two problems.

\subsection{Existence of a flat network output}
\label{sec:level1}

Our results on network flatness rest on the fact that~\eqref{eq:network} decouples along the eigenvectors of $L$, exactly for linear agents and, at synchronized equilibria, to first order for nonlinear ones.

\begin{lemma}[Network flatness for linear agents]\label{lem:linear_net}
Let $a(x_i)=\alpha^\top x_i, i \in \mathcal V$ and $b\in \R_{>0}$. Then $\Sg$ is flat if and only if $(L,\e_\ell)$ is controllable.
\end{lemma}

\begin{proof}
For $a(x_i)=\alpha^\top x_i$ and $b\in \R_{>0}$, \eqref{eq:network} is the linear system $\dot x=\big(I_N\otimes M-b L\otimes \e_q\e_1^\top\big)x+b\,(\e_\ell\otimes \e_q)u_G$ with $M:=A_c+\e_q\alpha^\top$. Since $L$ is symmetric, it admits an orthonormal diagonalization $L=Q\Lambda Q^\top$. In the coordinates $\tilde x=(Q^\top\otimes I_q)x$, the system decouples into the $N$ subsystems
\[
\dot{\tilde x}_i=M_i\tilde x_i+b[Q]_{\ell i}\,\e_q\,u_G,\quad M_i:=M-b\lambda_i \e_q\e_1^\top,
\]
for $i \in \mathcal V$.
Each $M_i$ is in companion form with characteristic polynomial $p_i(s)=s^q-\alpha^\top(1,s,\dots,s^{q-1})^\top+b\lambda_i$, so $(M_i,\e_q)$ is controllable and every eigenvalue of $M_i$ has geometric multiplicity one. Moreover, $p_i(s)-p_j(s)=b(\lambda_i-\lambda_j)$ is constant, so the spectra of $M_i$ and $M_j$ intersect if and only if $\lambda_i=\lambda_j$ because any common root has to satisfy $p_i(s)-p_j(s)=0$. Using the Hautus criterion~\cite{Hautus1969}, the decoupled system is therefore controllable if and only if all $\lambda_i$ are distinct and $[Q]_{\ell i}\neq0$ for all $i\in \mathcal V$. Since $[Q]_{\ell i}=\e_\ell^\top q_i$ for the $i$-th column $q_i$ of $Q$, an eigenvector of $L$, these are exactly the conditions of Lemma~\ref{lem:leaderfollower}, so the decoupled system is controllable if and only if $(L,\e_\ell)$ is. The assertion follows from Lemma~\ref{lem:lin_flat}. Related results on diffusively coupled linear agents are given in~\cite{Zhang2014}.
\end{proof}

\begin{proposition}[Necessary condition]\label{prop:necessary}
{Let Assumptions~\ref{ass:agent} and \ref{ass:graph} hold}. If $\Sg$ is locally flat around a synchronized equilibrium $x^\ast_c$, $c\in\mathbb S$, then $(L,\e_\ell)$ is controllable.
\end{proposition}

\begin{proof}
Write~\eqref{eq:network} as $\dot x=F(x)+G(x)u_G$ as in \eqref{eq:network_affine}. At $x^\ast_c$, every agent has the state $\xi^\ast:=c\,\e_1$, so that $L(I_N\otimes \e_1^\top)x^\ast_c=c\,L\one_N=0$. Hence the derivative of $B(x)=\diag(b(x_1),\dots,b(x_N))$ does not contribute to the Jacobian of the coupling term at $x^\ast_c$, and
\[
\begin{aligned}
\frac{\partial F}{\partial x}(x^\ast_c)&=I_N\otimes(A_c+\e_q\alpha^\top)-\beta L\otimes \e_q\e_1^\top,\\
G(x^\ast_c)&=\beta(\e_\ell\otimes \e_q),
\end{aligned}
\]
with $\alpha:=\nabla a(\xi^\ast)$ and $\beta:=b(\xi^\ast)>0$ by Assumption~\ref{ass:agent}. This is the linear network of Lemma~\ref{lem:linear_net}. Since $c\in\mathbb S$, the state $x^\ast_c$ is an equilibrium for $u_G=0$, so by Lemma~\ref{lem:siso_flat}(a) and Remark~\ref{rem:lin_nec}, local flatness around $x^\ast_c$ implies that this linearization is controllable, i.e., flat by Lemma~\ref{lem:lin_flat}. By Lemma~\ref{lem:linear_net}, this holds only if $(L,\e_\ell)$ is controllable.
\end{proof}

Proposition~\ref{prop:necessary} allows the classical graph-theoretic obstructions to leader--follower controllability to be carried over directly to flatness. The most prominent of these are symmetries. An automorphism of $\G$ is a relabeling of the nodes that preserves the edges and their weights~\cite{Godsil2001}; algebraically, it is a permutation matrix $P$ with $PL=LP$, and it fixes the input node if $P\e_\ell=\e_\ell$. Any two nodes exchanged by an automorphism fixing $\ell$ are influenced identically by $u_G$, so their states cannot be steered independently. Rahmani et al.~\cite{Rahmani2009} showed that $(L,\e_\ell)$ is then uncontrollable; see also~\cite{Egerstedt2012}. Combined with the spectral condition recalled in Section~\ref{sec:prelim}, this yields the following consequence of Proposition~\ref{prop:necessary}, which is illustrated in Fig.~\ref{fig:graphs}(d).

\begin{corollary}[Symmetry obstruction]\label{cor:symmetry}
If $\G$ admits a nontrivial automorphism fixing $\ell$, or if $L$ has a repeated eigenvalue, then $\Sg$ is not locally flat around any synchronized equilibrium, regardless of $a$ and $b$.
\end{corollary}

For linear agents, controllability of $(L,\e_\ell)$ is necessary and sufficient. The following example shows that this is no longer true for nonlinear agents, already for the path $\mathcal P_4$ driven from an interior node.

\begin{example}[Controllability is not sufficient]\label{ex:counter}
Consider scalar agents $\dot y_i=(1+y_i^2)\,u_i$, i.e., $q=1$, $a\equiv0$ and $b(\xi)=1+\xi^2$, coupled over $\mathcal P_4$ with unit weights and input node $\ell=2$. Then $(L,\e_2)$ is controllable, and $x^\ast_c$ is a synchronized equilibrium for every $c\in\R$. Since $\Sg$ has dimension four, local flatness requires by Lemma~\ref{lem:siso_flat}(a) that $\Delta_2$ be involutive. A direct symbolic computation gives
\begin{multline*}
\det\big[\,G,\ \mathrm{ad}_FG,\ \mathrm{ad}_F^2G,\ [\mathrm{ad}_FG,\mathrm{ad}_F^2G]\,\big]\\
=2\,(x_1-x_3)\,(1+x_1^2)^2(1+x_2^2)^5(1+x_3^2)^2(1+x_4^2),
\end{multline*}
which is nonzero whenever $x_1\neq x_3$, in particular at points arbitrarily close to any $x^\ast_c$. At such points $[\mathrm{ad}_FG,\mathrm{ad}_F^2G]\notin\Delta_2$, while $\dim\Delta_2=3$ near $x^\ast_c$ by controllability of the linearization. Thus $\Delta_2$ is not involutive near $x^\ast_c$, and $\Sg$ is not locally flat around any synchronized equilibrium, although $(L,\e_2)$ is controllable. For general weights and $b$, the factor $2(x_1-x_3)=b'(x_1)-b'(x_3)$ becomes $w_{12}\,b'(x_1)-w_{23}\,b'(x_3)$, so the same conclusion holds whenever $b$ is not affine near $c$. Since all vector fields involved are polynomial in $x$ and in the derivatives of $a$ and $b$, a failure of involutivity at a single point for one choice of agents implies failure for generic agents.
\end{example}

\subsection{Flat task outputs and graph distance}

We now turn to Problem~\ref{prob:node}. By Lemma~\ref{lem:siso_flat}(b), the task output is a flat output exactly when its relative degree equals the state dimension $Nq$, so it suffices to compute this relative degree. Intuitively, the input has to propagate from node $\ell$ to node $m$ through the graph, and each agent along the way contributes $q$ integrators.

\begin{proposition}[Relative degree of the task output]\label{prop:reldeg}
Let Assumptions~\ref{ass:agent} and~\ref{ass:graph} hold and let $d=\dist(\ell,m)$. Then the task output $y_G$ of  \eqref{eq:network_affine} has uniform relative degree
\[
\rho = q\,(d+1)
\]
with respect to $u_G$ on all of $\R^{Nq}$, and
\begin{equation}\label{eq:hfg}
\mathcal L_G\mathcal L_F^{\rho-1}h(x)=\sum_{\pi\in\Pi_{\ell m}}\;\prod_{\{u,v\}\in\pi}w_{uv}\;\prod_{j\in\pi}b(x_j)\;>\;0,
\end{equation}
where $\Pi_{\ell m}$ is the set of shortest paths from $\ell$ to $m$ in $\G$.
\end{proposition}

\begin{proof}
To investigate the relative degree of the network output, the proof proceeds in three steps: We first revisit \eqref{eq:network_affine}, then we analyze the algebraic structure of its drift vector field, before drawing upon the graph properties to show the assertion.

First write~\eqref{eq:network} as in \eqref{eq:network_affine}: $\dot x=F(x)+G(x)u_G$ with output $h(x)=x_{m,1}$. Assign to each state $x_{i,k}$ the \emph{level} $\nu(i,k):=q\,\dist(i,m)+k-1$.
By~\eqref{eq:agent} and~\eqref{eq:coupling}, the components of the drift vector field $F(x)$ in \eqref{eq:network_affine} are
\[
\begin{aligned}
F_{i,k}&=x_{i,k+1},\qquad k=1,\dots,q-1,\\
F_{i,q}&=a(x_i)+b(x_i)\sum_{j\in\N_i}w_{ij}\,(x_{j,1}-x_{i,1}).
\end{aligned}
\]
Intuitively, the level $\nu(i,k)$ counts how often the output $y_G$ has to be differentiated with respect to time until $x_{i,k}$ appears.

For $\mu\ge-1$, let $\mathcal C_\mu$ denote the set of smooth functions that depend only on states of level at most $\mu$ (so $\mathcal C_{-1}$ consists of the constant functions).
Since $|\dist(j,m)-\dist(i,m)|\le1$ for neighbors $i,j$, the component $F_{i,k}$ belongs to $\mathcal C_{\nu(i,k)+1}$ and is affine in its states of level exactly $\nu(i,k)+1$. The reason is that the coupling between agents is linear, as seen in the term $F_{i,q}$ above. For $k <q$ these are $x_{i,k+1}$, with coefficient $1$, and, if $k=q$, the states $x_{j,1}$ of the neighbors $j$ with $\dist(j,m)=\dist(i,m)+1$, with coefficients $w_{ij}\,b(x_i)\in\mathcal C_{\nu(i,q)}$. In particular, Lemma~\ref{lem:Inclusion} in the Appendix shows that $\mathcal L_F$ maps $\mathcal C_\mu$ into $\mathcal C_{\mu+1}$, i.e., $\mathcal L_F \mathcal C_\mu\subseteq \mathcal C_{\mu+1}$.

For a node $i$, let $\Pi_{im}$ be the set of shortest paths from $i$ to $m$ and
\[
\kappa_i(x):=\sum_{\pi\in\Pi_{im}}\;\prod_{\{u,v\}\in\pi}w_{uv}\prod_{j\in\pi\setminus\{i\}}b(x_j).
\]
The function $\kappa_i$ depends only on agents $v$ with $\dist(v,m)<\dist(i,m)$, whose states have level at most $q\,\dist(i,m)-1$; hence $\kappa_i\in\mathcal C_{s-1}$ whenever $\nu(i,k)=s$ for some $k$. We show by induction that, for all $s\ge0$,
\begin{equation}\label{eq:Lfs}
\mathcal L_F^sh=\sum_{(i,k):\,\nu(i,k)=s}\kappa_i\,x_{i,k}+R_s,\qquad R_s\in\mathcal C_{s-1}.
\end{equation}
For $s=0$, this holds with $\kappa_m=1$ and $R_0=0$. Suppose that~\eqref{eq:Lfs} holds for some $s$. Then
\[
\mathcal L_F^{s+1}h=\sum_{\nu(i,k)=s}\big((\mathcal L_F\kappa_i)\,x_{i,k}+\kappa_i\,F_{i,k}\big)+\mathcal L_FR_s .
\]
Since $\kappa_i,R_s\in\mathcal C_{s-1}$, the terms $(\mathcal L_F\kappa_i)\,x_{i,k}$ and $\mathcal L_FR_s$ belong to $\mathcal C_s$. By the structure of $F$, the states of level $s+1$ enter only through $\kappa_iF_{i,k}$, namely as $\kappa_i\,x_{i,k+1}$ if $k<q$, and as $\kappa_i\,w_{ij}\,b(x_i)\,x_{j,1}$ for the neighbors $j$ of $i$ with $\dist(j,m)=\dist(i,m)+1$ if $k=q$. Collecting terms, the coefficient of $x_{j,1}$ is $\sum_i w_{ij}\,b(x_i)\,\kappa_i$, summed over the neighbors $i$ of $j$ with $\dist(i,m)=\dist(j,m)-1$. This equals $\kappa_j$, because every shortest path from $j$ to $m$ consists of an edge to such a neighbor followed by a shortest path from that neighbor to $m$. All remaining terms belong to $\mathcal C_s$, which proves~\eqref{eq:Lfs} for $s+1$.

Finally, observe that $G(x)$ has only one non-zero entry and thus $\mathcal L_G\varphi=b(x_\ell)\,\frac{\partial\varphi}{\partial x_{\ell,q}}$ for any smooth function $\varphi$, and $\nu(\ell,q)=q(d+1)-1=\rho-1$. Since $\mathcal L_F^sh\in\mathcal C_s$ by~\eqref{eq:Lfs}, we obtain $\mathcal L_G\mathcal L_F^sh=b(x_\ell)\,\frac{\partial \mathcal L_F^sh}{\partial x_{\ell,q}}=0$ for all $s<\rho-1$ and all $x$, while $\mathcal L_G\mathcal L_F^{\rho-1}h=b(x_\ell)\,\kappa_\ell(x)$, which is~\eqref{eq:hfg}. This expression is positive on $\R^{Nq}$, since $\Pi_{\ell m}\neq\emptyset$ by Assumption~\ref{ass:graph}, $w_{ij}>0$, and $b>0$ by Assumption~\ref{ass:agent}.
\end{proof}

For $q=1$, $a\equiv0$, $b\equiv1$ and unit weights, Proposition~\ref{prop:reldeg} recovers~\cite{Briegel2011}: the relative degree is $\dist(\ell,m)+1$, and the high-frequency gain is the number of shortest $\ell$--$m$ paths.

\begin{lemma}[Flat task outputs]\label{lem:path}
Let Assumptions~\ref{ass:agent} and~\ref{ass:graph} hold. Then $y_G=y_m$ is a flat output of $\Sg$ if and only if $\G$ is the path graph $\mathcal P_N$ and $\ell$ and $m$ are its two end nodes. In this case the flatness maps $\Phi,\Psi$ are defined globally on $\R^{Nq}$.
\end{lemma}

\begin{proof}
By Lemma~\ref{lem:siso_flat}(b) and Proposition~\ref{prop:reldeg}, which yields a uniform relative degree on all of $\R^{Nq}$, the output $y_G$ is a flat output if and only if $q(\dist(\ell,m)+1)=Nq$, i.e., $\dist(\ell,m)=N-1$. In this case, a shortest path $\ell=v_0,v_1,\dots,v_{N-1}=m$ visits all $N$ nodes, and any further edge $\{v_i,v_j\}$ with $j>i+1$ would yield the shorter path $v_0,\dots,v_i,v_j,\dots,v_{N-1}$. Hence $\G=\mathcal P_N$ with end nodes $\ell$ and $m$; the converse is immediate.

On such a path, every node has a different distance to $m$, so the levels $\nu(i,k)$ from the proof of Proposition~\ref{prop:reldeg} take each value $0,\dots,Nq-1$ exactly once. By~\eqref{eq:Lfs}, the $s$-th component of $\Phi(x):=(h,\mathcal L_Fh,\dots,\mathcal L_F^{Nq-1}h)$ thus equals $\kappa_i\,x_{i,k}+R_s$, where $x_{i,k}$ is the unique state of level $s$, the functions $\kappa_i$ and $R_s$ depend only on states of lower level, and $\kappa_i>0$ on $\R^{Nq}$. These equations can be solved successively for the states, so $\Phi$ is a global diffeomorphism of $\R^{Nq}$. Since $y_G^{(s)}=\mathcal L_F^sh$ for $s<Nq$, this gives $\phi=\Phi^{-1}$ and $u_G=\big(y_G^{(Nq)}-\mathcal L_F^{Nq}h(x)\big)/\mathcal L_G\mathcal L_F^{Nq-1}h(x)$ with $x=\phi(\bar y_G^{[Nq-1]})$, which defines $\Psi$ globally because the denominator is positive by Proposition~\ref{prop:reldeg}.
\end{proof}

\begin{remark}[Flat inputs]\label{rem:flat_input}
Problem~\ref{prob:node} also has a dual reading: given the measured node $m$, how should the input be designed so that $y_m$ becomes a flat output? Since Lemma~\ref{lem:path} is symmetric in $\ell$ and $m$, a single-node input achieves this only on a path measured at an end node. If the input may act on several agents, however, the existence of such a \emph{flat input} is characterized by an observability rank condition alone, without any involutivity requirement~\cite{Waldherr08a}. This suggests that, around synchronized equilibria, a flat input for $y_m$ exists whenever $(L,\e_m)$ is controllable, for linear and nonlinear agents alike, at the price of an input that acts on many agents simultaneously. We leave a detailed analysis to future work.
\end{remark}

Note that the results of Section~\ref{sec:level1} do not involve the output node $m$, whereas Lemma~\ref{lem:path} shows that the task output is flat only in the very special case of a path measured at the end opposite the input. When network flatness holds but the task output is not flat, the flat output of the network need not be related to any particular agent, even for linear agents, as the following example shows.

\begin{example}[A hidden flat output]\label{ex:hidden}
Let $q=1$, $a\equiv0$, $b\equiv1$, and let $\G=\mathcal P_4$ with unit weights and input node $\ell=2$. The pair $(L,\e_2)$ is controllable, so $\Sg$ is flat by Lemma~\ref{lem:linear_net}. By Lemma~\ref{lem:lin_flat}, its flat output is, up to scaling,
\[
z=y_1-y_3-y_4,
\]
whose first three derivatives are independent of $u_G$, while $z^{(4)}$ depends on $u_G$ with coefficient $-1$. In contrast, the node outputs $y_1,\dots,y_4$ have relative degrees $2,1,2,3<4$ (Proposition~\ref{prop:reldeg}), so \emph{no} node output is a flat output, whichever node is measured. Nevertheless, trajectories of the whole network can be planned through $z$, even though $z$ is not available from any single measurement.
\end{example}

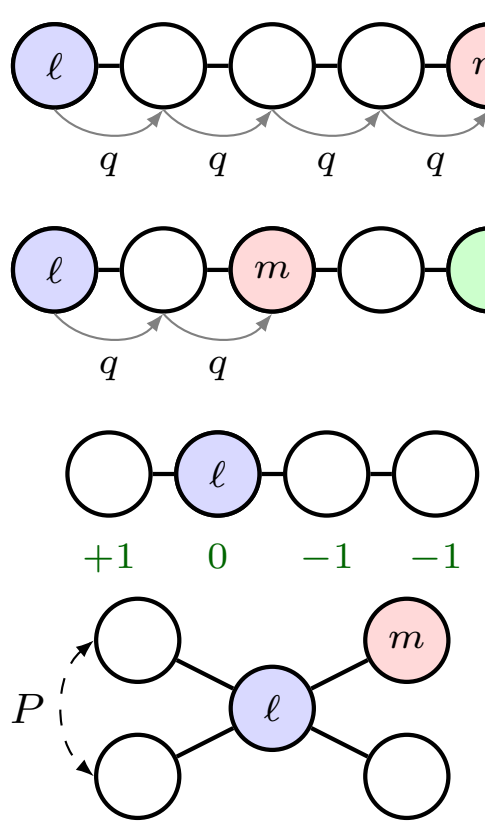
\begin{figure}[t]
\centering
\begin{tikzpicture}[font=\footnotesize]
\foreach \i in {1,...,5} {\node[agent] (a\i) at ({0.72*(\i-1)},0) {};}
\node[inn]  at (a1) {$\ell$};
\node[outn] at (a5) {$m$};
\draw[thick] (a1)--(a2)--(a3)--(a4)--(a5);
\foreach \i/\j in {1/2,2/3,3/4,4/5}
  {\draw[-{Latex[length=1.2mm]},gray] (a\i.south) to[bend right=50] node[below,font=\scriptsize,text=black]{$q$} (a\j.south);}
\node[anchor=west,text width=4.0cm] at (3.35,-0.1)
  {(a) Flat task output: $\rho=Nq$, the measured $y_G$ is the flat output (Lemma~\ref{lem:path}).};
\begin{scope}[yshift=-1.35cm]
\foreach \i in {1,...,5} {\node[agent] (b\i) at ({0.72*(\i-1)},0) {};}
\node[inn]  at (b1) {$\ell$};
\node[outn] at (b3) {$m$};
\node[flatn] at (b5) {};
\draw[thick] (b1)--(b2)--(b3)--(b4)--(b5);
\foreach \i/\j in {1/2,2/3}
  {\draw[-{Latex[length=1.2mm]},gray] (b\i.south) to[bend right=50] node[below,font=\scriptsize,text=black]{$q$} (b\j.south);}
\node[anchor=west,text width=4.0cm] at (3.35,-0.1)
  {(b) Network flatness only: $\Sg$ is flat with flat output $y_5$, but the measured $y_G$ has $\rho=3q<Nq$.};
\end{scope}
\begin{scope}[yshift=-2.7cm]
\foreach \i in {1,...,4} {\node[agent] (c\i) at ({0.36+0.72*(\i-1)},0) {};}
\node[inn] at (c2) {$\ell$};
\draw[thick] (c1)--(c2)--(c3)--(c4);
\foreach \i/\w in {1/{+1},2/{0},3/{-1},4/{-1}}
  {\node[below=0.5mm of c\i,font=\scriptsize,text=green!40!black] {$\w$};}
\node[anchor=west,text width=4.0cm] at (3.35,-0.1)
  {(c) Linear agents ($q=1$): network flatness only, flat output $z=y_1-y_3-y_4$ (weights in green; Ex.~\ref{ex:hidden}). Nonlinear agents: not flat (Ex.~\ref{ex:counter}).};
\end{scope}
\begin{scope}[yshift=-4.25cm]
\node[inn]   (d0) at (1.44,0)      {$\ell$};
\node[agent] (d1) at (0.55,0.45)  {};
\node[agent] (d2) at (0.55,-0.45) {};
\node[outn]  (d3) at (2.33,0.45)  {$m$};
\node[agent] (d4) at (2.33,-0.45) {};
\draw[thick] (d0)--(d1) (d0)--(d2) (d0)--(d3) (d0)--(d4);
\draw[{Latex[length=1.2mm]}-{Latex[length=1.2mm]},dashed] (d1.west) to[bend right=50] node[left]{$P$} (d2.west);
\node[anchor=west,text width=4.0cm] at (3.35,0)
  {(d) Neither: an automorphism $P$ fixes $\ell$, so $(L,\e_\ell)$ is uncontrollable and $\Sg$ is not flat (Cor.~\ref{cor:symmetry}).};
\end{scope}
\end{tikzpicture}
\caption{Flatness properties of networks; input node $\ell$ (blue), measured node $m$ (red), flat-output node (green). Each hop from $\ell$ towards $m$ adds $q$ to the relative degree of $y_G$ (Prop.~\ref{prop:reldeg}).}
\label{fig:graphs}
\end{figure}

\begin{remark}[Hierarchy of flatness properties]\label{rem:hierarchy}
The results of this section combine into
\[
\underbrace{y_G\text{ flat output}}_{\text{Definition \ref{def:flat_output}}}\;\Longrightarrow\;\underbrace{\Sg\text{ flat}}_{\text{Definition \ref{def:flat_network}}}\;\Longrightarrow\;(L,\e_\ell)\text{ controllable},
\]
where the second implication holds at synchronized equilibria (Proposition~\ref{prop:necessary}) and is an equivalence for linear agents (Lemma~\ref{lem:linear_net}) but not for nonlinear agents (Example~\ref{ex:counter}). The first implication is strict, see Fig.~\ref{fig:graphs}: on a path driven from an end node, the network is flat, but an interior task output is not (Fig.~\ref{fig:graphs}(b)), and for linear agents the flat output need not be the output of any agent (Example~\ref{ex:hidden}, Fig.~\ref{fig:graphs}(c)). In particular, collocated input and output ($\ell=m$) never yields a flat output for $N\ge2$.
\end{remark}

\section{Numerical Examples}
\label{sec:example}

We illustrate how network flatness enables trajectory planning through a single input. We consider second-order agents ($q=2$) with
\begin{equation}\label{eq:example_agent}
a(x_i)=-\eta_1x_{i,2}-\eta_2(x_{i,2})^3,\qquad b(x_i)=1+\varepsilon\,(x_{i,1})^2,
\end{equation}
$\eta_1=0.5$, $\eta_2=0.2$, $\varepsilon=0.5$, coupled over $\mathcal P_N$ with unit weights, $\ell=1$ and $y_G=y_N$, which is a flat output by Lemma~\ref{lem:path}; here $\mathbb S=\R$. On a path graph, the flat parametrizations $\Phi$ and $\Psi$ follow from a backward recursion: $y_N$ and its derivatives determine $x_N$, the dynamics of agent~$N$ can be solved for $y_{N-1}$, and so on towards the input node, each agent consuming $q$ further derivatives of $y_G$. We evaluate this recursion numerically by propagating the derivatives of the reference with Leibniz's rule. To transfer the network from $x^\ast_0$ to $x^\ast_1$ in time $T$, we plan
\begin{equation}\label{eq:yref}
y_{\rm ref}(t)=\sigma(t/T),\qquad \sigma(\tau)=I_\tau(Nq+1,Nq+1),
\end{equation}
where $I_\tau$ is the regularized incomplete beta function, extended by $0$ for $\tau<0$ and by $1$ for $\tau>1$. All derivatives of $y_{\rm ref}$ up to order $Nq$ vanish at $t=0$ and $t=T$, so $u_G=\Psi(\bar y_{\rm ref}^{[Nq]})$ vanishes outside $(0,T)$ and brings the network to rest at $x^\ast_1$.

\begin{figure}[t]
\centering
\centering
\includegraphics[width=\linewidth]{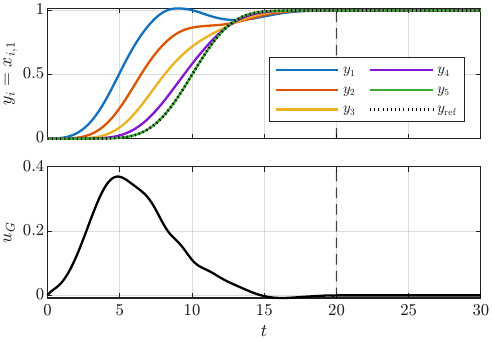}
\caption{Transition on $\mathcal P_5$ with input at node~1: outputs $y_i$, reference for $y_G=y_5$ (dotted), and $u_G$.}
\label{fig:ex1}
\end{figure}

Fig.~\ref{fig:ex1} shows the open-loop transition for $N=5$ ($Nq=10$) and $T=20$. The input acts on agent~1, which moves first and slightly overshoots, while agent~5 follows $y_{\rm ref}$. The network comes to rest at $x^\ast_1$ at $t=T$, and the simulated state deviates from the planned trajectory by less than $5\cdot10^{-10}$, which confirms the flatness maps numerically.

\begin{figure}[t]
\centering
\includegraphics[width=\columnwidth]{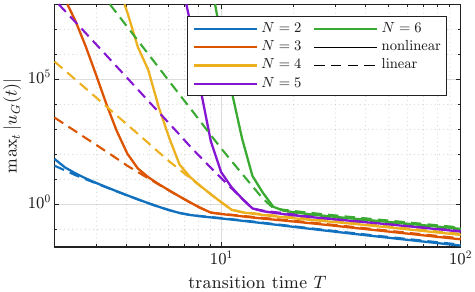}
\caption{Peak input for the transition $x^\ast_0\to x^\ast_1$ on $\mathcal P_N$ versus $T$. Solid: nonlinear agents~\eqref{eq:example_agent}; dashed: linear agents ($\eta_2=\varepsilon=0$).}
\label{fig:ex2}
\end{figure}

Fig.~\ref{fig:ex2} shows the peak input versus $T$. For linear agents ($\eta_2=\varepsilon=0$), $u_G$ is dominated by $y_{\rm ref}^{(Nq)}\propto T^{-Nq}$, and the fitted slopes $-4.07,\dots,-12.19$ agree with $-Nq$ for $N=2,\dots,6$. For nonlinear agents, the required input grows much faster below a critical transition time that increases with $N$, since each step of the recursion feeds the cubic drag of one agent into the next. Steering long nonlinear chains through one node is thus practical only over long horizons.

\begin{figure}[t]
\centering
\includegraphics[width=\linewidth]{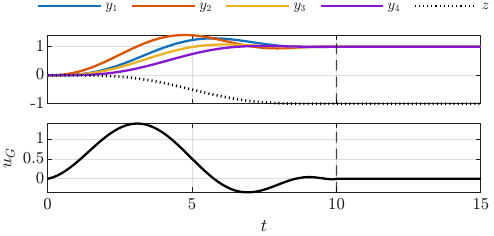}
\caption{Transition through the hidden flat output $z$ of Example~\ref{ex:hidden} (linear agents, $\mathcal P_4$, $\ell=2$).}
\label{fig:ex3}
\end{figure}

Finally, Fig.~\ref{fig:ex3} shows a transition through the hidden flat output of Example~\ref{ex:hidden}. Since $z=-c$ at $x^\ast_c$, we plan $z_{\rm ref}=-\sigma(t/T)$ with $Nq=4$ and $T=10$, and compute $u_G$ via Lemma~\ref{lem:lin_flat}. All agents reach consensus at $1$ at $t=T$ (error $4\cdot10^{-10}$), although no agent's output could have been used to plan the transition.

\section{Conclusions}
\label{sec:conclusion}

We have shown that flatness of a consensus network driven from a single node requires controllability of the Laplacian leader--follower pair, which is also sufficient for linear agents but not for nonlinear ones, whereas the task output is a flat output only if the network is a path with the input and the task output at its two ends.
An open problem is a graph-theoretic characterization of network flatness for nonlinear agents. Symbolic computations on small graphs suggest that, for generic nonlinear agents, the gap between network flatness and flat task outputs closes, i.e., that only paths driven from an end node are flat; for scalar agents, however, $\mathcal P_3$ driven from its center is flat whenever $w_{12}\neq w_{23}$. Further directions are trajectory planning under input constraints~\cite{Faulwasser2014}, multiple input and output nodes, directed graphs, nonlinear diffusive couplings, and the relation of Problem~\ref{prob:node} to \emph{flat inputs}~\cite{Waldherr08a}.

\section*{Appendix}
\begin{lemma}[$\mathcal L_F\mathcal C_\mu\subseteq\mathcal C_{\mu+1}$]\label{lem:Inclusion} ~\\
Consider \eqref{eq:network_affine} and \[T(F_{i,k}) = \left\{\begin{array}{l l}
    x_{i,k+1} & \quad \text{for } k <q \\
    \displaystyle \sum_{j\in\mathcal J_i^+} w_{ij}\,b(x_i)\,x_{j,1}\;+\;r_i(x) &
    \quad \text{for } k =q
\end{array}\right.,\]
 with index set  \[\mathcal J_i^+ := \{j\in\mathcal N_i\mid \dist(j,m)=\dist(i,m)+1\},\] and $r_i\in\mathcal C_{\nu(i,q)}$ reads
 \[
 \begin{aligned}
     r_i(x)&:=a(x_i)\;-\;b(x_i)\,[L]_{ii}\,x_{i,1}\;+\;b(x_i)\sum_{j \in \mathcal J_i^-} w_{ij}\,x_{j,1},\\
  &\quad  [L]_{ii}=\sum_{j\in\mathcal N_i}w_{ij}\\
 &\quad  \mathcal J_i^- := \{j\in\mathcal N_i \mid \dist(j,m)\le\dist(i,m)\}.
 \end{aligned}
 \]
Then the following statements hold:
\begin{itemize}
    \item[(i)] For all $\mu\ge-1$, the element-wise application of $\mathcal L_F$ to $\mathcal C_\mu$  gives $\mathcal L_F\mathcal C_\mu\subseteq\mathcal C_{\mu+1}$.
     \item[(ii)] For $\varphi\in\mathcal C_\mu$, we have
\begin{equation}\label{eq:LFtop}
  \mathcal L_F\varphi=\sum_{(i,k):\,\nu(i,k)=\mu}\frac{\partial\varphi}{\partial x_{i,k}}\,T(F_{i,k})+\tilde r(x),\ \ \tilde r\in\mathcal C_\mu.
\end{equation}
\item[(iii)] In particular, $\mathcal L_F\varphi$ is affine in the states of level exactly $\mu+1$, with coefficients in $\mathcal C_\mu$.
\end{itemize}
\end{lemma}

\begin{proof}
Part (i):
     Let $\varphi\in\mathcal C_\mu$. By the chain rule and since $\partial\varphi/\partial x_{i,k}=0$ whenever $\nu(i,k)>\mu$,
\begin{equation*}
  \mathcal L_F\varphi=\sum_{(i,k):\,\nu(i,k)\le\mu}\frac{\partial\varphi}{\partial x_{i,k}}\,F_{i,k}.
\end{equation*}
Each factor $\partial\varphi/\partial x_{i,k}$ depends only on states of level at most $\mu$ and therefore lies in $\mathcal C_\mu\subseteq\mathcal C_{\mu+1}$. Moreover, $F_{i,k}\in\mathcal C_{\nu(i,k)+1}\subseteq\mathcal C_{\mu+1}$ because $\nu(i,k)\le\mu$. Since $\mathcal C_{\mu+1}$ is closed under products and finite sums, $\mathcal L_F\varphi\in\mathcal C_{\mu+1}$.

Parts (ii) + (iii): For the refinement \eqref{eq:LFtop}, split the sum into the terms with $\nu(i,k)<\mu$ and those with $\nu(i,k)=\mu$. In the first group $F_{i,k}\in\mathcal C_{\nu(i,k)+1}\subseteq\mathcal C_\mu$ and $\partial\varphi/\partial x_{i,k}\in\mathcal C_\mu$, so the whole group lies in $\mathcal C_\mu$. In the second group, $F_{i,k}=T(F_{i,k})$ with $r_i\in\mathcal C_{\nu(i,q)}$ by the structure of~\eqref{eq:network_affine}, which is affine in level-$(\mu+1)$ states with coefficients in $\mathcal C_\mu$, plus a remainder in $\mathcal C_\mu$. Multiplying by $\partial\varphi/\partial x_{i,k}\in\mathcal C_\mu$ preserves both properties. Collecting all remainders in $\tilde r\in\mathcal C_\mu$ gives \eqref{eq:LFtop}.
\end{proof}

\section*{Acknowledgements}
During the preparation of this work, the authors used Claude (Anthropic; models Fable 5.1 and Opus 5.5) to structure the proofs of the results and to assist with editing the manuscript. After using this tool/service, the authors reviewed, restructured, and edited the content as needed and take full responsibility for the content of the publication.
\bibliographystyle{IEEEtran}
\bibliography{references}

\end{document}